\font\sets=msbm10.
\font\stampatello=cmcsc10.
\font\script=eusm10.

\def\1{{\bf 1}}
\def\oneQ{{\1_{_{[1,Q]}}}}
\def\oneR{{\1_{_{[1,R]}}}}
\def\oneH{{\1_{_{[-H,H]}}}} 
\def\gQ{{g_{_Q}}}
\def\muR{{\mu_{_R}}}
\def\wH{{w_{_H}}}
\def\uH{{u_{_H}}}

\def\sgn{{\rm sgn}}

\def\avesum{\sum_{x\sim N}}
\def\starsum{\mathop{\enspace{\sum}^{\ast}}}
\def\square{\hbox{\vrule\vbox{\hrule\phantom{s}\hrule}\vrule}}
\def\defineq{\buildrel{def}\over{=}}

\def\doublesum{\mathop{\sum\sum}}
\def\C{\hbox{\sets C}}
\def\N{\hbox{\sets N}}
\def\R{\hbox{\sets R}}
\def\Z{\hbox{\sets Z}}
\def\Corr{\hbox{\script C}}
\def\EssBdd{\ll_{\varepsilon}N^{\varepsilon}\,}
\def\Res{\mathop{{\rm Res}\,}}

\par
\centerline{\bf Sieve functions in arithmetic bands}
\bigskip
\par
\centerline{\stampatello g. coppola and m. laporta}
\bigskip
\par
\noindent {\bf Abstract}. An arithmetic function $f$ is called a {\it sieve function of range} $Q$, if its Eratosthenes transform $g=f\ast\mu$ is supported in $[1,Q]\cap\N$, where $g(q)\ll_{\varepsilon} q^{\varepsilon}$ ($\forall\varepsilon>0$). Here, we study the distribution of $f$ over short {\it arithmetic bands} $\cup_{1\le a\le H}\{n\in(N,2N]: n\equiv a\, (\bmod\,q)\}$, with $H=o(N)$, and give applications to both the correlations and to the so-called weighted Selberg integrals of $f$, on which we have concentrated our recent research.

\bigskip

\par
\centerline{\bf 1. Introduction and statement of the results.}
\smallskip
\par
\noindent
An arithmetic function $f:\N \rightarrow \C$ is called a {\it sieve function} of {\it range} $Q$, if 
$$
f(n)=\sum_{{d|n}\atop {d\le Q}}g(d), 
$$
where $g:\N \rightarrow \C$ is {\it essentially bounded}, namely $g(d)\ll_{\varepsilon} d^{\varepsilon}$, $\forall\varepsilon>0$. As usual, $\ll$ is Vinogradov's notation, synonymous to Landau's $O$-notation. In particular,  $\ll_{\varepsilon}$ means that the implicit constant might depend on an arbitrarily small and positive real number $\varepsilon$, which might change at each occurrence. 
When $f$ is the convolution product of $g$ and 
the constantly $1$ function, i.e. 
$$
f(n)=(g\ast \1)(n)=\sum_{d|n}g(d),
$$
\par
\noindent
we say, with Wintner [W], that $g$ is the {\it Eratosthenes transform} of $f$.
Observe that $f=g\ast \1$ is a sieve function whenever it is assumed that
$g$ is essentially bounded and vanishes outside $[1,Q]$ for some $Q\in\N$, that is to say, the Eratosthenes transform of $f$ is the restriction $\gQ\defineq g\cdot \oneQ$ (hereafter, $\1_B$ denotes the indicator function of the set $B\cap\Z$).
Moreover, by the M\"{o}bius inversion formula it turns out that $f=g\ast \1$ is essentially bounded 
if and only if so is $g$.  

Sieve functions are ubiquitous in analytic number theory. 
For example, the truncated divisor sum $\Lambda_R$, exploited by
Goldston in [G], is a linear combination of sieve functions of range $R$ (see Sect.4). Compare also [C2] for more examples of sieve functions. 
However, the reader is cautioned that by a sieve function some authors simply mean any sieve-related function that often arises  within the theory of sieve methods (see [DH]).

The first author has intensively investigated symmetry properties of sieve functions in short intervals through the study of their {\it correlations} and the associated {\it Selberg integrals} ([C1], [C2] and [CL1]). Here we wish to relate such a study to the distribution of a sieve function in modular arithmetic {\it short bands}. More precisely, for given positive integers $q,N,H$ 
we search for non-trivial bounds on the {\it total} ({\it balanced}) {\it value} of $f$ in {\it arithmetic bands} modulo $q$ defined as
$$
T_f(q,N,H)\defineq\sum_{a\le H}\sum_{{n\sim N}\atop {n\equiv a\, (\!\!\bmod q)}}f(n)
-{H\over q}\sum_{n\sim N}f(n), 
$$
\par
\noindent
where $n\sim N$ means that $n\in(N,2N]\cap \N$ (hereafter, we omit $a\ge 1$ in sums like $\sum_{a\le H}$). In particular, given any $N,H\in\N$, we prove that (see the remark after Theorem 1) for every real sieve function $f$ of range $Q\ll N$ and every $q\ll N$ one has 
$$
T_f(q,N,H)\ll_{\varepsilon}N^{\varepsilon}(N/q+q+Q). 
\leqno{(1)}
$$
\par
\noindent
It transpires from our method that similar bounds can be immediately established for {\it weighted} versions of the above problem, namely
$$
T_{w,f}(q,N,H)\defineq\sum_{0\le |a|\le H}w(a)\sum_{{n\sim N}\atop {n\equiv a\, (\!\!\bmod q)}}f(n)-{1\over q}\sum_{0\le |h|\le H}w(h)\sum_{n\sim N}f(n),
$$
\par				
\noindent
whenever $w:\R\rightarrow \R$ is a piecewise-constant {\it weight}. Indeed, it is plain that $T_f(q,N,H)=T_{u,f}(q,N,H)$ involves the {\it unit step} weight
$$
u(h)\defineq
\cases{1 & if $h>0$\cr 
0 & otherwise.\cr 
}
$$
\par
\noindent
However, we give more general conditions on $w$ to treat $T_{w,f}(q,N,H)$. First, let us set 
$$
\wH(h)\defineq w\cdot \oneH(h)=
\cases{
w(h) & if $h\in[-H,H]\cap\Z$\cr 
0 & otherwise,\cr 
}
$$
$$
{\cal L}^{1}_\ell(\widehat{\wH})\defineq {1\over\ell}\sum_{{j<\ell}\atop{(j,\ell)=1}}\Big|\widehat{\wH}\Big({j\over\ell}\Big)\Big|,
\quad \hbox{where}\quad 
\widehat{\wH}(\beta)\defineq \sum_{0\le |h|\le H}w(h)e(h\beta),
$$
\par
\noindent
(hereafter, $e(\alpha)\defineq e^{2\pi i\alpha}\ \forall \alpha\in \R$, and $(j,\ell)\defineq{\rm g.c.d.}(j,\ell)$, as usual in number theory). Thus, we can write
$$
\sum_{a}\wH(a)\sum_{{n\sim N}\atop {n\equiv a\, (\!\!\bmod q)}}f(n)=
{\widehat{\wH}(0)\over q}\sum_{n\sim N}f(n)+T_{w,f}(q,N,H) 
$$
and state our first result.
\smallskip
\par
\noindent {\stampatello Theorem 1}. {\it Let $q,N,H,Q$ be positive integers such that $q\ll N$ and $Q\ll N$, as $N\to \infty$. For every sieve function $f:\N \rightarrow \R$ of range $Q$ and every weight $w:\R \rightarrow \R$ one has
$$
T_{w,f}(q,N,H)\ll_{\varepsilon} N^{\varepsilon}\Big({N\over q}+q+Q\Big)\max_{{\ell>1}\atop{\ell|q}}{\cal L}^{1}_\ell(\widehat{\wH}).
$$
}
\smallskip
\par
\noindent {\stampatello Remark 1}. By taking $w=u$ and recalling $\Vert r\Vert\defineq \min_{n\in \Z}|r-n|$, $\forall r\in \R$, we have $\forall \ell>1$ [Da, Ch.25], 
$$
{\cal L}^{1}_\ell(\widehat{\uH})= 
{1\over\ell}\sum_{{j<\ell}\atop{(j,\ell)=1}}\Big|\sum_{h\le H}e\Big(h{j\over\ell}\Big)\Big|\ll
{1\over\ell}\sum_{{j<\ell}\atop{(j,\ell)=1}}{1\over {\Vert j/\ell\Vert}}\ll \sum_{j\le \ell/2}{1\over j}\ll \log \ell.
$$
Therefore, 
$(1)$ follows immediately from Theorem 1. 
\par
Another remarkable instance concerns the  {\it correlation} of $\wH$ given by
$$
W_H(a)\defineq\doublesum_{{h_1 \quad h_2}\atop {h_2-h_1=a}}\wH(h_1)\wH(h_2)
=\sum_{{0\le |h|\le H}\atop {0\le |h-a|\le H}}w(h)w(h-a).
$$
\par
\noindent
Note that $W_H$ vanishes outside \enspace $[-2H,2H]$. Moreover, uniformly in $\beta \in [0,1]$, 
$$
\widehat{W_H}(\beta)=\sum_{0\le |h|\le 2H}W_H(h)e(h\beta) 
=\sum_h \doublesum_{m-n=h}\wH(m)\wH(n)e(h\beta)
=\Big|\sum_r \wH(r)e(r\beta)\Big|^2=|\widehat{\wH}(\beta)|^2.
$$
\par
\noindent
Besides revealing that not all the weights are correlations of other weights, this yields
$$
\widehat{W_H}(0)=
\widehat{\wH}(0)^2\ll H^2,
$$
\par
\noindent
when $\wH$ is uniformly bounded as $H\to\infty$. Moreover, if $\wH$ also satisfies the inequality
$$
{\cal L}^{2}_\ell(\widehat{\wH})\defineq 
{1\over {\ell^2}}\sum_{j<\ell}\Big| \widehat{\wH}\left({j\over \ell}\right)\Big|^2\ll  {H\over \ell},\ \forall\ell\geq 1, 
\leqno{(2)}
$$
\par				
\noindent
then 
$$
{\cal L}^{1}_\ell(\widehat{W_H})={1\over\ell}\sum_{{j<\ell}\atop{(j,\ell)=1}}\widehat{W_H}\Big({j\over\ell}\Big)\leq \ell{\cal L}^{2}_\ell(\widehat{\wH})\ll H,\ \forall\ell\geq 1.
$$
\hfill $\diamond$
\par
\noindent
According to [CL2], a uniformly bounded weight $\wH$ (as $H\to\infty$) is said to be {\it good}, if it satisfies 
$(2)$. Thus, the following result is immediately established in a completely analogous way to
the proof of Theorem 1. 
\smallskip
\par
\noindent {\stampatello Corollary 1}. {\it Let $q,N,H,Q$ be positive integers such that $q\ll N$ and $Q\ll N$, as $N\to \infty$. For every sieve function $f:\N \rightarrow \R$ of range $Q$ and every good weight $w:\R \rightarrow \R$ one has
$$
\sum_{a}W_H(a)\sum_{{n\sim N}\atop {n\equiv a\, (\!\!\bmod q)}}f(n)
={{\widehat{W_H}(0)}\over q}\sum_{n\sim N}f(n)+O_{\varepsilon}\Big(N^{\varepsilon}H\Big({N\over q}+q+Q\Big)\Big),
$$
\par
\noindent
where $W_H$ is the correlation of $\wH$.
}
\smallskip
\par
\noindent {\stampatello Remark 2}. Though analogous definitions can be easily formulated for a complex weight $w$ (with the only exception of $W_H$, whose definition has to be modified by taking the complex conjugate of 
$\wH(h_2)$), here we stick to real weights and real sieve functions for simplicity.
\hfill $\diamond$
\par
\noindent {\stampatello Remark 3}. From [CL2] (see Propositions 2 and 3 there) it turns out that, beyond the unit step function $u$ defined above, other remarkable examples of good weights are the {\it sign} function 
and the {\it Cesaro} weight, respectively defined as
$$
\sgn(h)\defineq
\cases{
0 & if $h=0$\cr 
h/|h| & otherwise,\cr
}
\enspace 
C_H(h)\defineq
\cases{
1-|h|/H & if $|h|\le H$\cr 
0 & otherwise.\cr
}
$$
\par
\noindent
Since 
$$
C_H(h)={1\over H} \sum_{t\le H-|h|}1 ={1\over H} \doublesum_{{m,n\le H}\atop {m-n=h}}1,
$$ 
\par
\noindent
then $HC_H$ is the correlation of $\uH$, and consequently $\widehat{C_H}(0)=\widehat{\uH}(0)^2/H=H$. We conclude that Corollary 1 is non-trivial for 
$\wH=\uH$, yielding
$$
\sum_{a}C_H(a)\sum_{{n\sim N}\atop {n\equiv a\, (\!\!\bmod q)}}f(n)
={H\over q}\sum_{n\sim N}f(n)+O_{\varepsilon}\Big(N^{\varepsilon}\Big({N\over q}+q+Q\Big)\Big).
$$
\hfill $\diamond$
\par
\noindent {\stampatello Remark 4}. The main terms in the formul\ae\ furnished by Theorem 1 and Corollary 1 can be explicitly related to the 
Eratosthenes transform of $f=\gQ\ast\1$, with $Q\ll N$. Indeed, 
$$
\sum_{{n\sim N}\atop {n\equiv a\, (\!\!\bmod q)}}f(n)=\sum_{{n\sim N}\atop {n\equiv a\, (\!\!\bmod q)}}\sum_{d|n}\gQ(d)
=\sum_{d\le Q}g(d)\sum_{{{n\sim N}\atop {n\equiv a\, (\!\!\bmod q)}}\atop n\equiv 0\, (\!\!\bmod d)}1
=\sum_{d\le Q\atop (d,q)|a}g(d)\sum_{{n\sim N/d}\atop {nd\equiv a\, (\!\!\bmod q)}}1= 
$$
$$
=\sum_{d\le Q\atop (d,q)|a}g(d)\left( {N\over {dq}}(d,q)+O(1)\right)
={N\over q}\sum_{d\le Q\atop (d,q)|a}{g(d)\over d}(d,q)
 +O_{\varepsilon}\left(Q^{1+\varepsilon}\right). 
$$
\par
\noindent
In particular, for the {\it long intervals} we get the formula
$$
\sum_{n\sim N}f(n)=R_1(f)N+O_{\varepsilon}\left(Q^{1+\varepsilon}\right),
\leqno{(3)}
$$
\par				
\noindent
where the so-called first {\it Ramanujan coefficient} $R_1(f)$ is the mean value of $f$ (see Sect.1):
$$
R_1(f)\defineq\sum_{d\le Q}{{g(d)}\over d}=\lim_{x\to \infty}\Big(\sum_{d\le Q}{{g(d)}\over d}+{1\over x}\sum_{d\le Q}\!O\left(|g(d)|\right)\Big)
=\lim_{x\to \infty}{1\over x}\sum_{n\le x}f(n). 
$$
\par				
\noindent
On the other side, by taking $F$ as the Dirichlet series generating $f$, one has 
$$
\lim_{x\to \infty}{1\over x}\sum_{n\le x}f(n)=\Res_{s=1}F(s){{x^{s-1}}\over s}. 
$$
\par
\noindent
Since $f=\gQ\ast\1$ is a sieve function, then $F$ can be expressed in terms of 
the Riemann zeta function $\zeta$ and 
the Dirichlet polynomial generating its Eratosthenes transform, namely
$$
F(s)\defineq \sum_{n=1}^{\infty}{{f(n)}\over {n^s}}=\zeta(s)\sum_{d\le Q}{{g(d)}\over {d^s}}.
$$ 
\par				
\noindent
Note that the zeta function forces $F$ to have a simple pole at $s=1$, provided the $g$ series does not vanish at $s=1$. 
Thus, if $f=\gQ\ast\1$ is gauged by a weight $w$ in the {\it short interval} $[x-H,x+H]$ (i.e. $H=o(N)$, as $N\to \infty$), then it is natural to take the expected {\it mean value} of $\wH(n-x)f(n)$ for $N<x\le 2N$ to be (compare [CL]) 
$$
\widehat{\wH}(0)R_1(f)=\sum_a \wH(a)\sum_{d\le Q}{{g(d)}\over d}\quad (\hbox{that is independent of}\ x).
$$
\par
\noindent
Indeed, a basic tool for the study of the distribution of the sieve function $f$ in short intervals is its {\it weighted Selberg integral}
$$
J_{w,f}(N,H)\defineq \avesum \Big| \sum_n \wH(n-x)f(n)-\widehat{\wH}(0)R_1(f)\Big|^2, 
$$
whose non-trivial bounds might lead to results on the distribution of $f$ in {\it almost all} short intervals $[x-H,x+H]$, i.e. with $o(N)$ possible exceptions $x\in (N,2N]\cap\N$. Observe that the trivial bound for $J_{w,f}(N,H)$ is 
$N^{1+\varepsilon}H^2$, because $f$ is essentially bounded.
In [CL] and [CL2] we have investigated and exploited the link between $J_{w,f}(N,H)$ and the {\it correlation}
$$
\Corr_f(a)\defineq \sum_{n\sim N}f(n)f(n-a),
$$
in order to get non-trivial bounds under suitable conditions on $f$ and a good weight $w$.
\hfill $\diamond$
\medskip
\par
As a consequence of Theorem 1, we obtain a further result on such a link with a slight generalization. Let us define the correlation of real arithmetic functions $f_1$ and $f_2$ as
$$
\Corr_{f_1,f_2}(a)\defineq \sum_{n\sim N}f_1(n)f_2(n-a). 
$$
In such a context, we might refer to $\Corr_f=\Corr_{f,f}$ as the {\it autocorrelation} of $f$.
Since the {\it shift} $a$ is confined to $a\ll H$, 
the conditions $n\sim N$ and  $H=o(N)$ clearly yield $\max(n,n-a)\le 2N+|a|\le 3N$. Moreover,
if $f_1$ and $f_2$ are essentially bounded, then trivially $\Corr_{f_1,f_2}(0)\ll N^{1+\varepsilon}$, and for any 
$a\ll H$ one has 
$$
\Corr_{f_1,f_2}(a)=\doublesum_{{n_1\sim N\enspace n_2\sim N}\atop {n_2-n_1=a}}f_1(n_1)f_2(n_2)+O_{\varepsilon}\left(N^{\varepsilon}H\right) 
$$
\par
\noindent
(to be compared to the previous definition of the correlation of a weight).
\par
\noindent 
Correspondingly, the {\it mixed} weighted Selberg integral associated to the pair $(f_1,f_2)$ is 
(compare [C3])
$$
J_{w,(f_1,f_2)}(N,H)\defineq \avesum\prod_{j=1,2} \Big( \sum_n \wH(n-x)f_j(n)-\widehat{\wH}(0)R_1(f_j)\Big)\, .
$$
\par
\noindent
By applying Theorem 1  we obtain a {\it first generation} formula (consistently with the terminology of [CL]) for the correlation of the sieve functions $f_1$ and $f_2$, together with an estimate of the mixed weighted Selberg integral when these functions are gauged by a good weight $w$. 
\smallskip
\par				
\noindent {\stampatello Corollary 2}. {\it Let $N,H,Q_1,Q_2$ be positive integers with $Q_1\le Q_2\ll N$, as $N\to \infty$. For any real and essentially bounded arithmetic functions $g_1$ and $g_2$ supported in $[1,Q_1]$  and $[1,Q_2]$, respectively, one has 
$$
\sum_{a\le H}\Corr_{f_1,f_2}(a)=R_1(f_1)R_1(f_2)NH+O_{\varepsilon}\big(N^{\varepsilon}(N+Q_2^2+Q_1H)\big), 
$$
\par
\noindent
where $f_j=g_j\ast \1$ for $j=1,2$. Furthermore, if $H=o(N)$, as $N\to \infty$, and $w:\R \rightarrow \R$ is a good weight, then 
$$
J_{w,(f_1,f_2)}(N,H)\EssBdd (NH+Q_2H^2+Q_2^2 H+H^3). 
$$
}
\medskip
\par
\noindent  {\stampatello Remark 5}. For every real sieve function $f$ of range $Q\ll N$, this corollary gives
$$
\sum_{a\le H}\Corr_f(a)=R_1^2(f)NH+O_{\varepsilon}\big(N^{\varepsilon}(N+Q^2+QH)\big),
$$
$$
J_{w,f}(N,H)\EssBdd (NH+QH^2+Q^2 H+H^3). 
$$
\par
\noindent
We stress that such a bound for the weighted Selberg integral has been already established  in Theorem 3 of [CL2]. In Sect.2 we propose a much 
simpler proof through the new approach of the {\it arithmetic bands} formul\ae\ provided by Theorem 1. 
\par
Furthermore, from such an approach we find an important relation between weighted Selberg integrals and  the {\it total (weighted) content} of a sieve function $f$ of range $Q\ll N$, namely (see Lemma 2 and the proof of Corollary 2)
$$
J_{w,f}(N,H)\EssBdd \sum_{q\le Q}\left|T_{W,f}(q,N,H)\right|+N^{\varepsilon}H^2(Q+H), 
\leqno{(4)}
$$
\par
\noindent
where for the correlation of $\wH$ we set 
$$
T_{W,f}(q,N,H)\defineq\sum_{a}W_H(a)\sum_{{n\sim N}\atop {n\equiv a\, (\!\!\bmod q)}}f(n)
-{{\widehat{W_H}(0)}\over q}\sum_{n\sim N}f(n).
$$
\hfill $\diamond$
\par
\noindent	
Beyond Corollary 1, more generally,  given an essentially bounded $f$, a non-trivial  bound like
$$
\sum_{q\le Q}\left|T_{W,f}(q,N,H)\right|\ll N^{1-\delta}H^2\ \hbox{for some real}\  \delta>0
$$
\par
\noindent
might yield a non-trivial  bound of the same type for $J_{w,f}(N,H)$ (but not necessarily with the same {\it gain} $N^{\delta}$) by means of 
$(4)$. Analogous considerations hold for mixed weighted Selberg integrals. Rather surprisingly, in spite of the fact that the presence of absolute values in the total content seems to prevent it from further possible cancellation, the next theorem makes it clear that there are non-trivial bounds for (weighted) Selberg integrals, involving a sieve function $f$, of range $Q\ll N^{1-\delta}$ for some $\delta>0$, if and only if there are non-trivial results on the distribution of $f$ in short arithmetic bands. 
\smallskip
\par
\noindent {\stampatello Theorem 2}. {\it Let $f:\N \rightarrow \R$ be a sieve function of range $Q\ll N^{1-\delta}$, for some $\delta>0$, and let $w:\R \rightarrow \R$ be such that $\wH$ is uniformly bounded for any $H\ll N^{1-\delta}$, as $N\to \infty$.
\par	
\noindent 
I) The following three assertions are equivalent: 
\item{i)} a non-trivial  bound holds for $\displaystyle{\sum_{q\le Q}\left|T_{W,f}(q,N,H)\right|}$ 
\smallskip
\item{ii)} a non-trivial  bound holds for $J_{w,f}(N,H)$
\smallskip
\item{iii)} a non-trivial  bound holds for $J_{w,(f,f_1)}(N,H)$, where $f_1$ is any sieve function of range $Q$. 
\smallskip
\par
\noindent 
II) If $N^{\delta/2}\ll H\ll N^{1-\delta}$, as $N\to \infty$, then the following assertions are equivalent: 
\item{iv)} a non-trivial  bound holds for $\displaystyle{\sum_{q\le Q}\left|T_{f}(q,N,H)\right|}$ 
\smallskip
\item{v)} a non-trivial  bound holds for the Selberg integral 
$$
J_{f}(N,H)\defineq \sum_{x\sim N}\Big|\sum_{x<n\le x+H}f(n)-R_1(f)H\Big|^2. 
$$
}
\par
\noindent
Note that in $iv)$ a non-trivial bound is meant to be of the type $N^{1-\delta}H$ for some $\delta>0$. 
\bigskip

After a short section on some further notation and basic formul\ae, in Sect.2 we 
give the necessary lemmata for our theorems and for Corollary 2, whose proofs 
constitute the fourth section, whereas we omit the proof of Corollary 1, it being
completely analogous to the proof of Theorem 1. 
In Sect.4 we specialize the results of the present article to the aforementioned function $\Lambda_R$.
The last section is devoted to a comparison between classical results in arithmetic progressions and ours in arithmetic bands. 

\bigskip

\par
\centerline{\bf 2. Further notation and standard properties}
\smallskip
\par
\noindent
As already mentioned, we omit $a\ge 1$ in sums like $\sum_{a\le X}$. For the same sake of brevity, at times we write $n\equiv a\, (q)$ in place of $n\equiv a\, (\bmod\, q)$. 
Thus, the well-known {\it orthogonality of additive characters},
\par 
\centerline {$e_q(r)\defineq e(r/q)=e^{2\pi ir/q}$,\ ($q\in \N$,  $r\in \Z$),} 
\par
\noindent
can be written as
$$
{1\over q}\sum_{j\, (q)}e_q(j(n-m))={1\over q}\sum_{j\le q}e_q(j(n-m))=
\cases{
1 & if $n\equiv m\ (\bmod \; q)$\cr 
0 & otherwise\cr 
}
$$
\par
\noindent
since the sum is over a complete set of residue classes $j\ (\bmod\, q)$. 
\par
\noindent
We write  \thinspace ${\displaystyle \starsum_{j(q)} }$ \thinspace to mean that the sum is over a complete set of reduced residue classes $(\bmod \; q)$, i.e. the set $\Z_q^*$ of \enspace $1\le j\le q$ such that $(j,q)=1$.
In particular, the {\it Ramanujan sum} is written as
$$
c_q(n)\defineq \starsum_{j(q)}e_q(jn).
$$
\par
\noindent
Without further references, we will appeal to the well-known inequality (see [Da, Ch.25])
$$
\sum_{V_1<v\le V_2}e(v\alpha)\ll \min\Big(V_2-V_1,{1\over\Vert\alpha\Vert}\Big).
$$
Recalling that \thinspace $\1(n)\defineq 1,\ \forall n\in \N$, we set
$$
\1_{d|n}\defineq
\cases{
1 & if $d|n$\cr 
0 & otherwise.\cr 
}
$$
\par
\noindent
Consequently, the aforementioned orthogonality of characters becomes
$$
\1_{d|n}={1\over d}\sum_{j'(d)}e_d(j'n)={1\over d}\sum_{\ell|d}\sum_{{j'(d)}\atop {(j',d)=d/\ell}}e_d(j'n)
={1\over d}\sum_{\ell|d}c_{\ell}(n). 
$$
\par
\noindent
Therefore, one has the following {\it Ramanujan expansion} of a sieve function $f=\gQ\ast\1$:
$$
f(n)=\sum_{d|n}\gQ(d)
=\sum_{d\le Q}g(d)\1_{d|n}
=\sum_{d\le Q}{{g(d)}\over d}\sum_{\ell|d}c_{\ell}(n)
=\sum_{\ell \le Q}\sum_{{d\le Q}\atop {d\equiv 0\, (\ell)}}{{g(d)}\over d}c_{\ell}(n)
=\sum_{\ell \le Q}R_{\ell}(f)c_{\ell}(n), 
$$
\par				
\noindent
where we have introduced the so-called $\ell-$th {\it Ramanujan coefficient} of $f$, i.e.
$$
R_{\ell}(f)\defineq \sum_{d\equiv 0\, (\ell)}{{\gQ(d)}\over d}.
$$
\par
\noindent
The hypothesis that
$g$ is essentially bounded yields the bound
$$
R_{\ell}(f)\ll {1\over {\ell}}\sum_{m\le {Q\over {\ell}}}{{|g(\ell m)|}\over m}
\ll_{\varepsilon} {{Q^{\varepsilon}}\over {\ell}}\sum_{m\le {Q\over {\ell}}}{1\over m}
\ll_{\varepsilon} {{Q^{\varepsilon}}\over {\ell}}. 
\leqno{(5)}
$$
\par
\noindent
We refer the reader to [ScSp] and [W] for more extensive accounts on the theory of the Ramanujan expansions. 

\bigskip

\par
\centerline{\bf 3. Lemmata}
\smallskip
\par
\noindent
Here we state and prove two lemmas that are interesting in their own right. The first lemma is required to prove Theorem 1, while the second one is invoked within the proofs of Corollary 2 and Theorem 2. To this end, analogously to the exponential sums for the weights already introduced in Sect.1, we set 
$$
\widehat{f}(\alpha)\defineq \sum_{n\sim N}f(n)e(n\alpha) 
\quad 
(\alpha \in \R).
$$
\par
\noindent
Notice that now we can write
$$
T_{w,f}(q,N,H)=\sum_{a}\wH(a)\sum_{{n\sim N}\atop {n\equiv a\, (q)}}f(n)-\widehat{\wH}(0){\widehat{f}(0)\over q},
$$
$$
T_{W,f}(q,N,H)=\sum_{a}W_H(a)\sum_{{n\sim N}\atop {n\equiv a\, (q)}}f(n)-\widehat{W_H}(0){\widehat{f}(0)\over q},
$$
\par
\noindent
while the formula 
$(3)$ becomes 
$$
\widehat{f}(0)=R_1(f)N+O_{\varepsilon}(Q^{1+\varepsilon}). 
$$
\par
\noindent
The first lemma gives a similar relation between the $\ell-$th Ramanujan coefficient of $f$ and $\widehat{f}(\alpha)$, when $\alpha=j/\ell$ is any non-integer rational  with $(j,\ell)=1$. Note that such a formula is not a straightforward consequence of Wintner's criterion (see VIII.2 of [ScSp]). 
\smallskip
\par
\noindent {\stampatello Lemma 1.} {\it Let $f$ be a sieve function of range $Q\ll N$, with $Q,N\in\N$. Then
$$
\widehat{f}(j/\ell)=R_{\ell}(f)N+O_{\varepsilon}((\ell Q)^{\varepsilon}(Q+\ell)), 
\quad \enspace \forall \ell>1, \forall j\in\Z_{\ell}^{*}. 
$$
}
\smallskip
\par
\noindent {\stampatello Proof.} By assuming that $f=\gQ\ast\1$
with an essentially bounded $g$, we write
$$
\widehat{f}(j/\ell)=\sum_{d}\gQ(d)\sum_{v\sim {N\over d}}e_{\ell}(jdv)
=\sum_{d\equiv 0\, (\ell)}\gQ(d)\Big({N\over d}+O(1)\Big)+O\Big(\sum_{d\not \equiv 0\, (\ell)}{{|\gQ(d)|}\over {\Vert jd/\ell\Vert}}\Big). 
$$
\par
\noindent
Since 
$$
\sum_{d\equiv 0\, (\ell)}\gQ(d)\Big({N\over d}+O(1)\Big)
=R_{\ell}(f)N + O_{\varepsilon}\Big(Q^{\varepsilon}\Big({Q\over {\ell}}+1\Big)\Big)\ ,
$$
\par
\noindent
then the lemma is proved whenever we show that 
$$
\sum_{d\le Q\atop{d\not \equiv 0\, (\ell)}}{{1}\over {\Vert jd/\ell\Vert}}\ll_\varepsilon\ell^{\varepsilon}(Q+\ell).
$$
\par				
\noindent
To this end, it suffices to observe that
$$
\sum_{d\le Q\atop{d\not \equiv 0\, (\ell)}}{{1}\over {\Vert jd/\ell\Vert}}\le 
\sum_{0<|r|\le \ell/2}\sum_{d\le Q\atop{jd\equiv r\, (\ell)}}
{{1}\over {|r/\ell|}}\ll_\varepsilon
\ell\sum_{r\le \ell/2}{{1}\over r}\Big({Q\over {\ell}}+1\Big).
$$
The proof is completed.\hfill $\square$ 
\smallskip
\par
\noindent  {\stampatello Remark 6}. Note that the formula of the above lemma is non-trivial when 
$\ell, Q\ll N^{1-\delta}$, for some $\delta>0$. Moreover, it is easy to see that it holds uniformly with respect to $j\in\Z_{\ell}^{*}$.	
\hfill $\diamond$
\smallskip
\par
Let us turn our attention to the next lemma. As already mentioned in Sect.1, by means of an elementary dispersion method, in [CL], Lemma 7, we established a link between weighted Selberg integrals and autocorrelations of an arithmetic function $f$ gauged by a weight $w$ such that $\wH$ is bounded, as $H\to \infty$. Under the further hypothesis that  the sieve function $f$ and the weight $w$ are real, the formula of the aforementioned lemma becomes 
$$
J_{w,f}(N,H)=\sum_{0\le |a|\ll H}W_H(a)\Corr_f(a)-2\widehat{\wH}(0)R_1(f)\sum_{n\le 3N}f(n)\avesum \wH(n-x)+\avesum\left|\widehat{\wH}(0)R_1(f)\right|^2 + 
$$
$$
+O_\varepsilon\left(H^3N^{\varepsilon}\right). 
$$
\par
\noindent
Similarly, for the mixed weighted Selberg integral of sieve functions $f_1,f_2$ we have
$$
J_{w,(f_1,f_2)}(N,H)=\sum_a W_H(a)\Corr_{f_1,f_2}(a)-\widehat{W_H}(0)R_1(f_1)R_1(f_2)N-\widehat{\wH}(0)\times 
\leqno{(6)}
$$
$$
\times \Big(R_1(f_1)\avesum \Delta_2(x)+R_1(f_2)\avesum \Delta_1(x)\Big)
 +O_\varepsilon\left(H^3N^{\varepsilon}\right), 
$$
\par
\noindent
where we set \enspace $\Delta_j(x)\defineq\sum_n \wH(n-x)f_j(n)-\widehat{\wH}(0)R_1(f_j)$. By using such a formula we prove the next lemma, where $J_{w,(f_1,f_2)}(N,H)$ is expressed in terms of arithmetic bands of $f_1$ or $f_2$. 
\smallskip
\par
\noindent {\stampatello Lemma 2.} {\it Let $g_1$ and $g_2$ be real and essentially bounded arithmetic functions supported in $[1,Q_1]$  and $[1,Q_2]$, respectively, with $Q_1,Q_2\in \N$ such that $Q_1\le Q_2\ll N$, as $N\to \infty$. If $w:\R \rightarrow \R$ is such that $\wH$ is uniformly bounded, as $H\to \infty$, then one has 
$$
J_{w,(f_1,f_2)}(N,H)=\sum_{q\le Q_1}g_1(q)T_{W,f_2}(q,N,H)+O_{\varepsilon}\!\left(N^{\varepsilon}H^2(Q_2+H)\right)= 
$$
$$
=\sum_{q\le Q_2}g_2(q)T_{W,f_1}(q,N,H)+O_{\varepsilon}\!\left(N^{\varepsilon}H^2(Q_2+H)\right), 
$$
\par
\noindent
where we set $f_j=g_j\ast \1$, and $W_H$ is the correlation of $\wH$. 
}
\smallskip
\par
\noindent {\stampatello Proof.} First, let us write
$$
\avesum\sum_n \wH(n-x)f_j(n)=\sum_{n\sim N} f_j(n)\sum_{n-H\le x\le n+H} w(n-x)+O_{\varepsilon}\left(N^{\varepsilon}H^2\right)
=\widehat{\wH}(0)\sum_{n\sim N} f_j(n)+O_{\varepsilon}\left(N^{\varepsilon}H^2\right).
$$
\par
\noindent
Then, by arguing as in 
$(3)$ and recalling that $R_1(f_j)\ll_\varepsilon Q_j^\varepsilon$, we get
$$
\avesum\Delta_j(x)=\widehat{\wH}(0)\Big(\sum_{n\sim N} f_j(n)-R_1(f_j)N\Big)+O_{\varepsilon}\left(N^{\varepsilon}H^2\right)
\ll_{\varepsilon} N^{\varepsilon}H(Q_j+H). 
$$
\par
\noindent
Since $W_H$ is even and $Q_1\le Q_2\ll N$, the above formula 
$(6)$ yields 
$$
J_{w,(f_1,f_2)}(N,H)=\sum_a W_H(a)\Corr_{f_1,f_2}(a)-\widehat{W_H}(0)R_1(f_1)R_1(f_2)N
 +O_{\varepsilon}\left(N^{\varepsilon}H^2(Q_2+H)\right)= 
$$
$$				
=\sum_a W_H(a)\Corr_{f_2,f_1}(a)-\widehat{W_H}(0)R_1(f_1)R_1(f_2)N
 +O_{\varepsilon}\left(N^{\varepsilon}H^2(Q_2+H)\right). 
$$
\par
\noindent
Thus, we can stick to the first equality, apply 
$(3)$ to $f_1$ and write 
$$
\sum_a W_H(a)\Corr_{f_1,f_2}(a)-\widehat{W_H}(0)R_1(f_1)R_1(f_2)N= 
$$
$$
=\sum_a W_H(a)\sum_{n\sim N}f_1(n)\sum_{{q|n-a}\atop {q\le Q_2}}g_2(q)
 -\widehat{W_H}(0)\sum_{n\sim N}f_1(n)\sum_{q\le Q_2}{{g_2(q)}\over q}
  +O_{\varepsilon}\left(Q_2^{1+\varepsilon}H^2\right)= 
$$
$$
=\sum_{q\le Q_2}g_2(q)\Big(\sum_a W_H(a)\sum_{{n\sim N}\atop {n\equiv a\, (q)}}f_1(n)-\widehat{W_H}(0){\widehat{f_1}(0)\over q}\Big) 
 +O_{\varepsilon}\left(Q_2^{1+\varepsilon}H^2\right). 
$$
\par
\noindent
The lemma is completely proved.\hfill $\square$ 

\bigskip

\par
\centerline{\bf 4. Proofs of Theorems 1, 2 and Corollary 2}
\smallskip
\par
\noindent
{\stampatello Proof of Theorem 1.} By the orthogonality of additive characters we get 
$$
T_{w,f}(q,N,H)={1\over q}\sum_{a}\wH(a)\sum_{n\sim N}f(n)\sum_{j'\le q}e_q(j'(a-n))-{{\widehat{\wH}(0)}\over q}\widehat{f}(0)
={1\over q}\sum_{j'<q}\sum_{a}\wH(a)e_q(j'a)\widehat{f}\left(-j'/q\right)= 
$$
$$
={1\over q}\sum_{\ell>1\atop \ell|q}\starsum_{j\, (\ell)}\widehat{f}\left(-j/\ell\right)\widehat{\wH}(j/\ell), 
$$
\par
\noindent
where we have set $\ell=q/(j',q)$. By applying Lemma 1 and 
$(5)$ we see that 
$$
T_{w,f}(q,N,H)\ll_{\varepsilon} {1\over q}\sum_{{\ell>1}\atop {\ell|q}}\Big(|R_{\ell}(f)|N
                                 +(\ell Q)^{\varepsilon}(Q+\ell)\Big)\!\!\starsum_{j\, (\ell)}\Big|\widehat{\wH}\Big({j\over\ell}\Big)\Big|
\ll_{\varepsilon} {{Q^{\varepsilon}}\over q}\sum_{{\ell>1}\atop {\ell|q}}\Big({N\over {\ell}}
                   +Q\ell^{\varepsilon}+\ell^{1+\varepsilon}\Big)\ell {\cal L}^{1}_\ell(\widehat{\wH})\ll_{\varepsilon} 
$$
$$
\ll_{\varepsilon} N^{\varepsilon}\Big({N\over q}+Q+q\Big)\max_{{\ell>1}\atop{\ell|q}}{\cal L}^{1}_\ell(\widehat{\wH}). 
$$
\par
\noindent
The theorem is completely proved. \hfill $\square$
\bigskip
\par
\noindent
{\stampatello Proof of Corollary 2.} As already noticed in the proof of Lemma 2, we can write
$$
\Corr_{f_1,f_2}(a)=\sum_{n\sim N}f_1(n)f_2(n-a)=\sum_{n\sim N}f_1(n)\sum_{{q|n-a}\atop {q\le Q_2}}g_2(q)
=\sum_{q\le Q_2}g_2(q)\sum_{{n\sim N}\atop {n\equiv a\, (q)}}f_1(n).
$$
\par
\noindent
Thus, the formula 
$(3)$ and Theorem 1 yield 
$$
\sum_{a\le H}\Corr_{f_1,f_2}(a)=\sum_{q\le Q_2}g_2(q)\Big({H\over q}\widehat{f}_1(0)+T_{f_1}(q,N,H)\Big)= 
$$
$$
=H\Big(\sum_{q\le Q_2}{{g_2(q)}\over q}\Big)\left(R_1(f_1)N+O_{\varepsilon}\left(Q_1^{1+\varepsilon}\right)\right)
 +O_{\varepsilon}\Big(N^{\varepsilon}\sum_{q\le Q_2}\Big({N\over q}+q+Q_1\Big)\Big)= 
$$
$$
=R_1(f_1)R_1(f_2)NH+O_{\varepsilon}\big(N^{\varepsilon}(N+Q_2^2+Q_1H)\big), 
$$
\par
\noindent
that is the first formula of Corollary 2. In order to prove the stated inequality 
for the mixed weighted Selberg integral, it is enough to observe that Lemma 2 and the hypothesis
$Q_1\le Q_2\ll N$ imply
$$
J_{w,(f_1,f_2)}(N,H)=\sum_{q\le Q_2}g_2(q)T_{W,f_1}(q,N,H)+O_{\varepsilon}\big(N^{\varepsilon}H^2(Q_2+H)\big)\ll_{\varepsilon} 
$$
$$				
\ll_{\varepsilon} N^{\varepsilon}\sum_{q\le Q_2}\left|T_{W,f_1}(q,N,H)\right|+N^{\varepsilon}H^2(Q_2+H).
$$
\par
\noindent
Whence the conclusion follows immediately from Corollary 1. \hfill $\square$ 

\bigskip

Before going to the proof of Theorem 2, let us remark explicitly that 
$(4)$ is plainly a particular case of the latter inequality. Moreover, it transpires from the previous proof that, for every real and essentially bounded arithmetic function $g$ supported in $[1,Q]$, with $Q\ll N$, one has
$$
\sum_{a\le H}\Corr_f(N)=R_1(f)^2 NH+\sum_{q\le Q}g(q)T_f(q,N,H)+O_{\varepsilon}\left(N^{\varepsilon}QH\right), 
\leqno{(7)}
$$
\par
\noindent
where we set $f=g\ast \1$.
\bigskip
\par
\noindent
{\stampatello Proof of Theorem 2.} For simplicity and without loss of generality, let us assume that, whatever the choice of an assertion among {\it i)-v)} as hypothesis, the gain of the non-trivial bound is always $N^{\delta}$.
\smallskip
\par
\noindent
Part I. $i) \Longrightarrow ii)$: as we said, let us suppose that 
$$
\sum_{q\le Q}\left|T_{W,f}(q,N,H)\right|\ll N^{1-\delta}H^2. 
$$
\par
\noindent
Thus, $ii)$ follows immediately from 
$(4)$, where $H^2(Q+H)\ll N^{1-\delta}H^2$ because of the hypotheses $H,Q\ll N^{1-\delta}$. 
\smallskip
\par
\noindent
$ii) \Longrightarrow iii)$: since we assume that $J_{w,f}(N,H)\ll N^{1-\delta}H^2$, then
by the Cauchy inequality and the trivial bound for $J_{w,f_1}(N,H)$ we get 
$$
J_{w,(f,f_1)}(N,H)\le\sqrt{J_{w,f}(N,H)}\sqrt{J_{w,f_1}(N,H)}
\EssBdd\sqrt{N^{1-\delta}H^2}\sqrt{NH^2}
\ll N^{1-\delta/3}H^2.
$$
\par
\noindent
$iii) \Longrightarrow i)$: after setting
$$
s_{W,f}(q)\defineq
\cases{
\sgn(T_{W,f}(q,N,H)) & if $1\le q\le Q$\cr 
0 & otherwise,\cr 
}
$$ 
\par
\noindent
it is readily seen that $f_1=s_{W,f}\ast \1$ is a sieve function of range $Q$. Thus, we can write 
$$
\sum_{q\le Q}\left|T_{W,f}(q,N,H)\right|=\sum_{q}s_{W,f}(q)T_{W,f}(q,N,H). 
$$
\par
\noindent
Now, by taking $g_1=s_{W,f}$ and $f_2=f$  in Lemma 2 we see that 
$$
\sum_{q\le Q}\left|T_{W,f}(q,N,H)\right|=J_{w,(f,f_1)}(N,H)+O_{\varepsilon}\left(N^{\varepsilon}H^2(Q+H)\right),
$$
\par
\noindent
where again $H^2(Q+H)$ is non-trivial. The first part of the theorem is completely proved. 
\smallskip
\par
\noindent
Part II. $iv) \Longrightarrow v)$: since $Q\ll N^{1-\delta}$ and we assume that
$$
\sum_{q\le Q}\left|T_{f}(q,N,H)\right|\ll N^{1-\delta}H, 
$$
\par
\noindent
then it is easily seen that the formula 
$(7)$ yields 
$$
\sum_{a\le t}\Corr_f(a)=R_1(f)^2 N[t]+O_{\varepsilon}\left(N^{1-\delta+\varepsilon}t\right)\quad\hbox{for all}\ 1\le t\le H,
$$
\par				
\noindent
where \enspace $[t]$ \enspace is the {\it integer part} of \enspace $t$. Thus, by partial summation we can write
$$
\sum_{1\le a\le H}(H-a)\Corr_f(a)=\int_{1}^{H}\sum_{a\le t}\Corr_f(a)dt
=\int_{1}^{H}\left(R_1(f)^2 N[t]+O_{\varepsilon}\left(N^{1-\delta+\varepsilon}t\right)\right)dt= 
$$
$$
={{R_1(f)^2}\over 2}NH^2+O_{\varepsilon}\left(N^{1+\varepsilon}H\right)+O_{\varepsilon}\left(N^{1-\delta+\varepsilon}H^2\right). 
$$
\par
\noindent
Now, since $\Corr_f(0)\ll_{\varepsilon} N^{1+\varepsilon}$, and for $1\leq a\leq H$  one has
$$
\Corr_f(-a)=\sum_{n\sim N}f(n)f(n+a)=\sum_{N+a<m\le 2N+a}f(m-a)f(m)
=\Corr_f(a)+O_{\varepsilon}\left(N^{\varepsilon}H\right),
$$
\par
\noindent
then 
$$
\sum_{0\le |a|\le H}(H-|a|)\Corr_f(a)=R_1(f)^2 NH^2+O\Big(N^{1-\delta/3}H^2\Big).
$$
\par
\noindent
By using this formula in 
$(6)$, where we take $W_H(a)=HC_H(a)=\max(H-|a|,0)$ (see Remark 3), we immediately obtain $J_f(N,H)\ll N^{1-\delta/3}H^2$. 
\smallskip
\par
\noindent
$v) \Longrightarrow iv)$: we suppose that $J_f(N,H)\ll N^{1-\delta}H^2$ and set
$$
s_{f}(q)\defineq
\cases{
\sgn(T_{f}(q,N,H)) & if $1\le q\le Q$\cr 
0 & otherwise,\cr
}
\quad\quad 
f_1\defineq s_f\ast \1.
$$ 
\par
\noindent
Thus, we can write
$$
\sum_{q\le Q}\left|T_f(q,N,H)\right|=\sum_{q}s_f(q)\Big(\sum_{a\le H}\sum_{{n\sim N}\atop {n\equiv a\, (\!\!\bmod q)}}f(n)
	                                                 -{H\over q}\sum_{n\sim N}f(n)\Big)= 
$$
$$
=\sum_{a\le H}\Big(\sum_{n\sim N}f(n)f_1(n-a)-R_1(f_1)\sum_{n\sim N}f(n)\Big)= 
$$
$$
=\sum_{a\le H}\sum_{N-a<x\le 2N-a}f(x+a)f_1(x)-R_1(f_1)R_1(f)NH
 +O_{\varepsilon}(N^{\varepsilon}QH)= 
$$
$$
=\sum_{x\sim N}f_1(x)\sum_{x<m\le x+H}f(m)-R_1(f_1)R_1(f)NH
 +O_{\varepsilon}(N^{\varepsilon}(Q+H)H)= 
$$
$$
=\sum_{x\sim N}f_1(x)\Big(\sum_{x<m\le x+H}f(m)-R_1(f)H\Big)
 +O_{\varepsilon}(N^{\varepsilon}(Q+H)H),
$$
\par
\noindent
where we have applied 
$(3)$ to both $f$ and $f_1$. Note that the $O$-term contribution is non-trivial because of hypotheses on $Q$ and $H$. In order to deal with the main term of the latter formula, after recalling that $f_1$ is essentially bounded, we apply the Cauchy inequality and the above assumption on $J_f(N,H)$ to get
$$
\sum_{x\sim N}f_1(x)\Big(\sum_{x<m\le x+H}f(m)-R_1(f)H\Big)\ll_{\varepsilon} N^{1/2+\varepsilon}\sqrt{J_f(N,H)}
\ll_{\varepsilon} N^{1+\varepsilon-\delta/2}H,
$$
\par
\noindent
which in turn yields
$$
\sum_{q\le Q}\left|T_f(q,N,H)\right|\ll N^{1-\delta/3}H.
$$
\par
\noindent
Theorem 2 is completely proved. \hfill $\square$ 

\vfill
\eject

\par				
\centerline{\bf 5. A remarkable truncated divisor sum}
\smallskip
\par
\noindent
Let us recall that the truncated divisor sum is defined in [G] as 
$$
\Lambda_R(n)\defineq \sum_{{d|n}\atop {d\le R}}\mu(d)\log(R/d)=
(\log R)\sum_{{d|n}\atop {d\le R}}\mu(d)-\sum_{{d|n}\atop {d\le R}}\mu(d)\log d, 
$$
\par
\noindent
so that $\Lambda_R$ is plainly a linear combination (with relatively {\it small} coefficients) of two sieve functions, whose Eratosthenes transforms are respectively the {\it restricted} M\"obius function, $\muR\defineq \mu\cdot\oneR$, and $\muR\cdot\log$. 
\par
After recalling also the well-known relations (see [Da])
$$
\sum_{d=1}^{\infty}{{\mu(d)\log d}\over d}=-1 
\quad 
\hbox{\rm and}
\quad
\sum_{d\le R}{{\mu(d)}\over d}, \sum_{d>R}{{\mu(d)\log d}\over d} \ll \exp\big(-c\sqrt{\log R}\big),
$$
\par
\noindent
(hereafter, $c>0$ is an unspecified constant), we see that
$$
R_1(\Lambda_R)=\sum_{d\le R}{{\mu(d)\log(R/d)}\over d}
=(\log R)\sum_{d\le R}{{\mu(d)}\over d}-\sum_{d\le R}{{\mu(d)\log d}\over d}
=1+O\big(\exp\big(-c\sqrt{\log R}\big)\big).
$$
\par
\noindent
Thus, the mean value formula 
$(3)$ gives 
$$
\sum_{n\sim N}\Lambda_R(n)=N+O\big(N\exp\big(-c\sqrt{\log R}\big)\big)+O_{\varepsilon}(N^{\varepsilon}R), 
$$
while, if $R\ll N$, a straightforward application of 
$(1)$ yields 
$$
\sum_{a\le H}\sum_{{n\sim N}\atop {n\equiv a\, (\!\!\bmod q)}}\Lambda_R(n)
={{NH}\over q}+O_{\varepsilon}\Big(N^{\varepsilon}\Big({N\over q}+q+R\Big)\Big)+O\big(N\exp\big(-c\sqrt{\log R}\big)\big).
$$
\noindent
In case the {\it level} $\lambda\defineq(\log R)/(\log N)$ is positive, i.e. $0<\lambda_0\le \lambda<1$ (for a fixed $\lambda_0$), we may replace $\log R$ by $\log N$ in the above formul\ae, where now $c=c(\lambda)$. Assuming that this is the case, Corollary 2 provides the following {\it first generation} formula for the correlation of $\Lambda_R$:
$$
\sum_{a\le H}\sum_{n\sim N}\Lambda_R(n)\Lambda_R(n-a)=NH+O\big(NH\exp\big(-c\sqrt{\log N}\big)\big)
                                                         +O_{\varepsilon}\big(N^{\varepsilon}(N+R^2+RH)\big). 
$$
\par
\noindent
It is worthwhile to remark that by following the classical approach in the literature the remainder term for the single correlation
is $\ll_{\varepsilon} N^{\varepsilon}R^2$, that trivially yields a remainder $\ll_{\varepsilon} N^{\varepsilon}R^2 H$ in the 
first generation formula above, whereas by our method we save $H$. 

\bigskip

\par
\centerline{\bf 6. Further comments}
\smallskip
\par
\noindent
The key of the present approach is that the correlation of a real sieve function $f=\gQ\ast\1$ can be written as
$$
\Corr_f(a)=\sum_{q\le Q}g(q)\sum_{{n\sim N}\atop {n\equiv a(\!\!\bmod \, q)}}f(n).
$$
In the literature (see [Ik], Ch.17), we find several studies of the distribution of 
an arithmetic function $f$ (not necessarily a sieve function) over primitive residue classes. Most results are focused on non-trivial bounds for the {\it error term}
$$
E_f(N;q,a)\defineq \sum_{{n\sim N}\atop {n\equiv a(\!\!\bmod \, q)}}f(n)-M_f(N;q,a)
$$ 
\par				
\noindent
for all $(a,q)=1$, provided $q$ is not too large. Here, $M_f(N;q,a)$ is
the expected {\it mean value} term.
Let us recall two major variants of the problem. The first one concerns the {\it Bombieri-Vinogradov} type mean
$$
\sum_{q\le Q}\max_{(a,q)=1}\left|E_f(N;q,a)\right|,
$$
\par
\noindent
for which we refer the reader to [M]. The second classical variant is
the {\it Barban-Davenport-Halberstam} type quadratic mean
$$
\sum_{q\le Q}\sum_{{a\le q}\atop {(a,q)=1}}E_f(N;q,a)^2.
$$
\par
\noindent
The latter has also a short interval version introduced by Hooley [Ho], that is
$$
\sum_{q\le Q}\sum_{{a\le \rho q}\atop {(a,q)=1}}E_f(N;q,a)^2,\ \hbox{where}\ \rho\to 0.
$$
\par
\noindent
In all such problems, the challenging issue is the {\it level} $\lambda \defineq (\log Q)/(\log N)$
of distribution of $f$ in arithmetic progressions (see [FI], \S 9.8 and \S22.1). For example, 
the celebrated Bombieri-Vinogradov Theorem gives a non-trivial bound for
$$
\sum_{q\le Q}\max_{(a,q)=1}\Big|\sum_{{n\sim N}\atop {n\equiv a(\!\!\bmod \, q)}}\Lambda(n)-{N\over {\varphi(q)}}\Big|,
\quad \hbox{\rm where} \enspace 
\varphi(q)\defineq |\{ a\le q, (a,q)=1\}|, 
$$
\par
\noindent
essentially with a level $\lambda=1/2$ (which seems to be a structural barrier at least for the distribution of primes). However, for many applications one can just deal with individual reduced class $a$ and take the sum over $q\le Q$, $(q,a)=1$. Indeed, by assuming that $a\not=0$, one can see that  it is possible to break the level $1/2$ for the {\it Bombieri-Friedlander-Iwaniec} type mean (see [FI], Theorem 22.1)
$$
\sum_{q\le Q\atop (a,q)=1}\Big|\sum_{{n\sim N}\atop {n\equiv a(\!\!\bmod \, q)}}\Lambda(n)-{N\over {\varphi(q)}}\Big|.
$$
\par
\noindent
Consistently with the present notation, the above formula for the correlation of a sieve function becomes
$$
\Corr_f(a)=\sum_{q\le Q}g(q)M_f(N;q,a)+\sum_{q\le Q}g(q)E_f(N;q,a),
$$ 
\par
\noindent
where, by recalling that $g(q)\ll_\varepsilon q^\varepsilon$, one has
$$
\sum_{q\le Q}g(q)E_f(N;q,a)\ll_\varepsilon Q^\varepsilon\sum_{q\le Q}\left|E_f(N;q,a)\right|.
$$
\par
\noindent
Thus, here for each individual residue $a$ we deal with a sum over $q\le Q$ without any further restriction. Then, it is not surprising that
a straight asymptotic
$$
\Corr_f(a)\sim \sum_{q\le Q}g(q)M_f(N;q,a) 
$$
\par
\noindent
has been proved for very few interesting instances of $f$, including the noteworthy case of the divisor function (see the third version of [CL] on {\tt arxiv} for a brief account on this matter). Better expectations for the {\it first generation} of correlation averages,
$$
\sum_{a\le H}\Corr_f(a),
$$
\par				
\noindent
are given substance by Corollary 2 (and by the alternative approach of Lemma 12 in [CL]). Furthermore, note that Theorem 2 concerns the average
$$
\sum_{q\le Q}\Big| \sum_{a\le H}E_f(N;q,a)\Big|,
$$
\par
\noindent
where, unlike the aforementioned means, the sums are taken over all the moduli $q\le Q$ and over a short interval of residue classes $a$, when $f$ is a sieve function of range $Q\ll N^{1-\delta}$ and $H\ll N^{1-\delta}$. The bound for the weighted Selberg integral given in Corollary 2 and its application through Theorem 2 allow $Q\ll \sqrt{NH}N^{-\varepsilon}$, that is to say, the level might go beyond $1/2$ when we deal with {\it not too short} intervals, e.g., $H\gg N^{3\varepsilon}$. 

\bigskip

\par
\noindent {\bf Acknowledgment.} This research started while the first author was a fellow \lq \lq Ing.Giorgio Schirillo\rq \rq \thinspace of the Istituto Nazionale di Alta Matematica (Italy).

\bigskip

\par
\centerline{\stampatello References}
\smallskip
\item{[C1]} G. Coppola - {\sl On the Correlations, Selberg integral and symmetry of sieve functions in short intervals} - J. Combinatorics and Number Theory {\bf 2.2}, Article 1,  91--105, 2010
\item{[C2]} G. Coppola - {\sl On the Correlations, Selberg integral and symmetry of sieve functions in short intervals}, II - Int. J. Pure Appl. Math. {\bf 58.3}, 281--298, 2010
\item{[C3]} G. Coppola - {\sl On some lower bounds of some symmetry integrals} - Afr. Mat. {\bf 25}, issue 1, 183--195, 2014
\item{[CL]} G. Coppola and M. Laporta - {\sl Generations of correlation averages} - Journal of Numbers Vol. {\bf 2014}, Article ID 140840, 1-13, 2014 (compare draft http://arxiv.org/abs/1205.1706) 
\item{[CL1]} G. Coppola and M. Laporta - {\sl On the Correlations, Selberg integral and symmetry of sieve functions in short intervals, III} -  (submitted), http://arxiv.org/abs/1003.0302
\item{[CL2]} G. Coppola and M. Laporta - {\sl Symmetry and short interval mean-squares} - (submitted), see draft online at http://arxiv.org/abs/1312.5701
\item{[Da]} H. Davenport - {\sl Multiplicative Number Theory} - 3rd edition, GTM {\bf 74}, Springer, New York, 2000
\item{[DH]} H.G. Diamond and H. Halberstam - {\sl A Higher-Dimensional Sieve Method  (With Procedures for Computing Sieve Functions by W.F. Galway)} -  Cambridge Tracts in Mathematics, Vol. {\bf 177}, Cambridge University Press, Cambridge, 2008
\item{[FI]}  J. Friedlander and H. Iwaniec - {\sl Opera de Cribro} - AMS Colloquium Publications, {\bf 57}, Providence, RI, 2010 
\item{[G]} D.A. Goldston - {\sl On Bombieri and Davenport's theorem concerning small gaps between primes} - Mathematika {\bf 39}, 10--17, 1992
\item{[Ho]} C. Hooley - {\sl On the Barban-Davenport-Halberstam Theorem. XI} - Acta Arith. {\bf 91}, no.{\bf 1}, 1--41, 1999
\item{[Ik]} H. Iwaniec and E. Kowalski - {\sl Analytic Number Theory} - AMS Colloquium Publications, {\bf 53}, Providence, RI, 2004 
\item{[M]} Y. Motohashi - {\sl An induction principle for the generalization of Bombieri's prime number theorem} - Proc. Japan Acad. {\bf 52}, no.{\bf 6},  273--275, 1976
\item{[ScSp]} W. Schwarz and J. Spilker - {\sl Arithmetical functions (An introduction to elementary and analytic properties of arithmetic functions and to some of their almost-periodic properties)} - London Mathematical Society Lecture Note Series, {\bf 184}, Cambridge University Press, Cambridge, 1994
\item{[W]} A. Wintner - {\sl Eratosthenian Averages} - Waverly Press, Baltimore, MD, 1943

\vfill
\eject

\par				
\noindent
{\stampatello Giovanni Coppola}
\par
\noindent
Postal address: Via Partenio 12, 
\par
\noindent
83100 Avellino (AV), ITALY
\par
\noindent
e-mail : giocop@interfree.it
\par
\noindent
wpage: www.giovannicoppola.name
\bigskip

\noindent
{\stampatello Maurizio Laporta}
\par
\noindent
Universit\`a degli Studi di Napoli "Federico II",
\par
\noindent
Dipartimento di Matematica e Applicazioni "R.Caccioppoli",
\par
\noindent
Complesso di Monte S.Angelo,
\par
\noindent
Via Cinthia, 80126 Napoli (NA), ITALY
\par
\noindent
e-mail : mlaporta@unina.it

\bye